\documentclass[12pt]{amsart}
\usepackage{geometry}
\usepackage{amsmath}
\usepackage{amsfonts}
\usepackage{amssymb}
\usepackage{mathrsfs}
\usepackage{eufrak}
\allowdisplaybreaks[2]

\newcommand{\equstart}{\begin{equation}\begin{aligned}}
\newcommand{\equend}{\end{aligned}\end{equation}}
\newcommand{\equstartu}{\begin{equation*}\begin{aligned}}
\newcommand{\equendu}{\end{aligned}\end{equation*}}
\newtheorem{theorem}{Theorem}
\newtheorem{lemma}{Lemma}
\newtheorem{corollary}{Corollary}
\newtheorem{proposition}{Proposition}
\theoremstyle{definition}

\newtheorem{Question}{Question}
\theoremstyle{remark}

\theoremstyle{plain}
\newtheorem{other}{Theorem}         
\newtheorem{otherc}[other]{Corollary}

\newenvironment{pf}{\noindent{\emph{Proof.}}}{$\Box$}
\newenvironment{Pf}{\noindent{\emph{Proof of}}}{$\Box$}

\newcommand{\hol}{{\mathcal Hol}}
\DeclareMathOperator{\og}{O}

\DeclareMathOperator{\op}{o} 
\newcommand{\ig}{\stackrel{\text{def}}{=}}
\def\D{{\mathbb D}}

\def\C{{\mathbb C}}

\def\Dp{{\mathcal D^p_{p-1}}}

\def\Dpa{{\mathcal D^p_{\alpha}}}

\begin{document}
\title[]
{Hankel matrices acting on the Hardy space $H^1$ and on Dirichlet
spaces}

\author[D.~Girela]{Daniel Girela}
\address{An\'alisis Matem\'atico,
Universidad de M\'alaga, Campus de Teatinos, 29071 M\'alaga, Spain}
\email{girela@uma.es}
\author[N.~Merch\'{a}n]{Noel Merch\'{a}n}
\address{An\'alisis Matem\'atico,
Universidad de M\'alaga, Campus de Teatinos, 29071 M\'alaga, Spain}
\email{noel@uma.es} 
\subjclass[2010]{Primary 47B35; Secondary 30H10.}
\date{}
\keywords{Hankel matrix, Generalized Hilbert operator, Hardy spaces,
Cauchy transforms, Weighted Bergman spaces, Dirichlet spaces,
Duality}

\begin{abstract} If $\,\mu \,$ is a finite positive Borel measure on the interval
$\,[0,1)$, we let $\,\mathcal H_\mu \,$ be the Hankel matrix $\,(\mu
_{n, k})_{n,k\ge 0}\,$ with entries $\,\mu _{n, k}=\mu _{n+k}$,
where, for $\,n\,=\,0, 1, 2, \dots $, $\mu_n\,$ denotes the moment
of order $\,n\,$ of $\,\mu $. This matrix induces formally the
operator $\,\mathcal{H}_\mu (f)(z)=
\sum_{n=0}^{\infty}\left(\sum_{k=0}^{\infty}
\mu_{n,k}{a_k}\right)z^n\,$ on the space of all analytic functions
$\,f(z)=\sum_{k=0}^\infty a_kz^k\,$, in the unit disc $\,\mathbb D
$. When $\,\mu \,$ is the Lebesgue measure on $\,[0,1)\,$ the
operator $\,\mathcal H_\mu\,$ is the classical Hilbert operator
$\,\mathcal H\,$ which is bounded on $\,H^p\,$ if $\,1<p<\infty $,
but not on $\,H^1$. J. Cima has recently proved that $\,\mathcal
H\,$ is an injective bounded operator from $\,H^1\,$ into the space
$\,\mathscr C\,$ of Cauchy transforms of measures on the unit
circle.
\par The operator $\,\mathcal H_\mu \,$ is known to be well defined on
$\,H^1\,$ if and only if $\,\mu \,$ is a Carleson measure and in
such a case we have that $\mathcal H_\mu (H^1)\subset \,\mathscr C$.
Furthermore, it is bounded from $\,H^1\,$ into itself if and only if
$\,\mu\,$ is a $1$-logarithmic $1$-Carleson measure.
\par In this paper we prove that when
 $\,\mu\,$ is
a $1$-logarithmic $1$-Carleson measure then $\,\mathcal H_\mu \,$
actually maps $\,H^1\,$ into the space of Dirichlet type $\,\mathcal
D^1_0\,$. We discuss also the range of $\,\mathcal H_\mu\,$ on
$\,H^1\,$ when $\,\mu \,$ is an $\alpha $-logarithmic $1$-Carleson
measure ($0<\alpha <1$). We study also the action of the operators
$\,\mathcal H_\mu \,$ on Bergman spaces and on Dirichlet spaces.
\end{abstract}
\thanks{This research is supported in part by a grant from \lq\lq El Ministerio de
Econom\'{\i}a y Competitividad\rq\rq , Spain (MTM2014-52865-P) and
by a grant from la Junta de Andaluc\'{\i}a (FQM-210). The second
author is also supported by a grant from \lq\lq El Ministerio de
Educaci\'{o}n, Cultura y Deporte\rq\rq , Spain (FPU2013/01478).}
\maketitle
\section{Introduction and main results}\label{intro}
 Let $\,\D=\{z\in\C:|z|<1\}\,$
denote the open unit disc in the complex plane $\,\C$, $\partial
\mathbb D\,$ will be the unit  circle. The space of all analytic
functions in $\,\D\,$ will be denoted by $\,\hol (\mathbb D)$. We
also let $\,H^p\,$ ($0<p\le \infty $) be the classical Hardy spaces.
We refer to \cite{D} for the notation and results regarding Hardy
spaces.
\par For $\,0<p<\infty \,$ and $\,\alpha >-1\,$
the weighted Bergman space $\,A^p_\alpha \,$ consists of those
$\,f\in \hol (\mathbb D)\,$ such that
$$\Vert f\Vert _{A^p_\alpha }\,\ig\, \left ((\alpha +1)\int_\mathbb D(1-\vert z\vert ^2)^{\alpha }\vert f
(z)\vert ^p\,dA(z)\right )^{1/p}\,<\,\infty .$$ Here, $\,dA\,$
stands for the area measure on $\,\mathbb D$, normalized so that the
total area of $\,\mathbb D\,$ is $\,1$. Thus
$\,dA(z)\,=\,\frac{1}{\pi }\,dx\,dy\,=\,\frac{1}{\pi
}\,r\,dr\,d\theta $. The unweighted Bergman space $\,A^p_0\,$ is
simply denoted by $\,A^p$. We refer to \cite{DS,HKZ,Zhu} for the
notation and results about Bergman spaces.\par The space of
Dirichlet type $\,\Dpa\,$ ($0<p<\infty \,$ and $\,\alpha
>-1$) consists of those $f\,\in \hol (\mathbb D)\,$ such that $\,f^\prime
\in A^p_\alpha $. In other words, a function $\,f\in \hol (\mathbb
D)\,$ belongs to $\,\Dpa\,$ if and only if
$$\Vert f\Vert _{\Dpa}\,\ig\,\vert f(0)\vert \,+\,\left ((\alpha +1)\int_{\mathbb D}(1-\vert z\vert ^2)^{\alpha }\vert f^\prime
(z)\vert ^p\,dA(z)\right )^{1/p}\,<\,\infty .$$
\par The Hilbert matrix is the infinite matrix $\,\mathcal H=\left
(\frac{1}{k+n+1}\right )_{k,n\ge 0}$. It induces formally an
operator, called the Hilbert operator, on spaces of analytic
functions as follows:
\par If $\,f\in \hol (\D )$, $f(z)=\sum_{n=0}^\infty a_nz^n$, then we
set \begin{equation}\label{def-H}\mathcal Hf(z)=\sum_{n=0}^\infty
\left (\sum_{k=0}^\infty \frac{a_k}{n+k+1}\right )z^n,\quad
z\in\mathbb D,\end{equation} whenever the right-hand side of
(\ref{def-H}) makes sense for all $\,z\in \D\,$ and the resulting
function is analytic in $\,\D$. We define also
\begin{equation}\label{def-I}\mathcal
If(z)=\int_0^1\frac{f(t)}{1-tz}\,dt,\quad z\in\mathbb
D,\end{equation} if the integrals in the right-hand side of
(\ref{def-I}) converge for all $\,z\in\mathbb D\,$ and the resulting
function $\,\mathcal If\,$ is analytic in $\,\mathbb D$. It is clear
that the correspondences $\,f\mapsto \mathcal Hf\,$ and $\,f\mapsto
\mathcal If\,$ are linear.
\par If $\,f\in H^1$, $f(z)=\sum_{n=0}^\infty a_nz^z$, then by the
Fej\'{e}r-Riesz inequality \cite[Theorem\,\@3.\,\@13, p.\,\@46]{D}
and Hardy's inequality \cite[p.\,\@48]{D}, we have
$$\int_0^1\vert f(t)\vert\,dt\le\pi \Vert f\Vert
_{H^1}\quad\text{and}\quad \sum_{n=0}^\infty \frac{a_n}{n+1}\le \pi
\Vert f\Vert _{H^1}.$$ This immediately yields that if $\,f\in
H^1\,$ then $\,\mathcal Hf\,$ and $\mathcal If\,$ are well defined
analytic functions in $\,\mathbb D\,$ and that, furthermore,
$\mathcal Hf=\mathcal If$.
\par Diamantopoulos and Siskakis \cite{DiS} proved that
$\,\mathcal H\,$ is a bounded operator from $\,H^p\,$ into itself if
$\,1<p<\infty $, but this is not true for $\,p=1$. In fact, they
proved that $\mathcal H\left (H^1\right )\nsubseteq H^1$. Cima
\cite{Cima} has recently proved the following result.
\begin{other}\label{th-cima} \begin{itemize}\item[(i)] The operator
$\,\mathcal H \,$ maps $\,H^1\,$ into the space $\,\mathscr C\,$ of
Cauchy transforms of measures on the unit circle $\,\partial \mathbb
D$. \item[(ii)] $\,\mathcal H:H^1\rightarrow \mathscr C\,$ is
injective. \end{itemize}\end{other}
\par We recall that if $\,\sigma \,$ is a finite complex Borel
measure on $\,\partial \mathbb D$, the Cauchy transform $\,C\sigma
\,$ is defined by
$$C\sigma (z)=\int_{\partial \mathbb D}\frac{d\sigma (\xi
)}{1-\overline {\xi }\,z},\quad z\in \mathbb D.$$ We let $\mathscr
M$ be the space of all finite complex Borel measure on $\,\partial
\mathbb D$. It is a Banach space with the total variation norm. The
space of Cauchy transforms is $\,\mathscr C=\{ C\sigma : \sigma \in
\mathscr M\}$. It is a Banach space with the norm $\,\Vert C\sigma
\Vert \ig \inf \{ \Vert \tau \Vert : C\tau =C\sigma \}$. We mention
\cite{CMR} as an excellent reference for the main results about
Cauchy transforms. We let $\,\mathcal A\,$ denote the disc algebra,
that is, the space of analytic functions in $\,\mathbb D\,$ with a
continuous extension to the closed unit disc, endowed with the
$\,\Vert \cdot \Vert_{H^\infty }$-norm. It turns out
\cite[Chapter\,\@4]{CMR} that $\,\mathcal A\,$ can be identified
with the pre-dual of $\,\mathscr C\,$ via the pairing
\begin{equation}\label{dual-A}\langle g, C\sigma \rangle\, \ig
\,\lim_{r\to 1}\frac{1}{2\pi }\int _0^{2\pi }g(re^{i\theta
})\overline {C\sigma (re^{i\theta })}\,d\theta .\end{equation} This
is the basic ingredient used by Cima to prove the inclusion
$\,\mathcal H(H^1)\subset \mathscr C$.
\par Now we turn to consider a class of operators which are
natural generalizations of the operators $\,\mathcal H\,$ and
$\,\mathcal I$. If $\,\mu\,$ is a finite positive Borel measure on
$\,[0, 1)\,$ and $\,n\, = 0, 1, 2, \dots $, we let $\,\mu_n\,$
denote the moment of order $\,n\,$ of $\,\mu\,\,$, that is, $\mu
_n=\int _{[0,1)}t^n\,d\mu (t),$ and we define $\,\mathcal H_\mu \,$
to be the Hankel matrix $\,(\mu _{n,k})_{n,k\ge 0}\,$ with entries
$\,\mu _{n,k}=\mu_{n+k}$. The measure $\,\mu\,$ induces formally the
operators $\,\mathcal I_\mu \,$ and $\,\mathcal H_\mu \,$ on spaces
of analytic functions as follows:
\[\mathcal I_\mu f(z)=\int_{[0,1)}\frac{f(t)}{1-tz}\,d\mu (t),\quad
\mathcal H_\mu f(z)=\sum_{n=0}^\infty \left (\sum_{k=0}^\infty
a_k\mu_{n+k}\right )z^n,\quad z\in\mathbb D,\] for
$\,f(z)=\sum_{n=0}^\infty a_nz^n\in \hol (\mathbb D)\,$ being such
that the terms on the right-hand sides make sense for all $\,z\in
\mathbb D$, and the resulting functions are analytic in $\,\mathbb
D$. If $\,\mu \,$ is the Lebesgue measure on $\,[0,1)\,$ the matrix
$\,\mathcal H_\mu \,$ reduces to the classical Hilbert matrix and
the operators $\,\mathcal H_\mu\,$ and $\,\mathcal I_\mu \,$ are
simply the operators $\,\mathcal H\,$ and $\,\mathcal I$.
\par If $\,I\subset \partial\D\,$ is an interval, $\,\vert I\vert \,$ will
denote the length of $\,I$. The \emph{Carleson square} $\,S(I)\,$ is
defined as $\,S(I)=\{re^{it}:\,e^{it}\in I,\quad 1-\frac{|I|}{2\pi
}\le r <1\}$.
\par If $\, s>0\,$ and $\,\mu\,$ is a positive Borel  measure on  $\,\D$,
we shall say that $\,\mu \,$
 is an $s$-Carleson measure
  if there exists a positive constant $\,C\,$ such that
\[
\mu\left(S(I)\right )\le C{|I|^s}, \quad\hbox{for any interval
$I\subset\partial\D $}.
\]
A $\,1\,$-Carleson measure will be simply called a Carleson measure.
We recall that Carleson \cite{Carl} proved that $\,H^p\,\subset
\,L^p(d\mu )\,$ ($0<p<\infty $) if and only if $\,\mu \,$ is a
Carleson measure (see also \cite[Chapter\,\@9]{D}).
\par
For \,$0\le \alpha <\infty \,$ and $\,0<s<\infty \,$ we say that a
positive Borel measure $\,\mu\,$ on $\,\D\,$ is an
 $\alpha$-logarithmic $s$-Carleson measure if there exists a positive
 constant $C$ such that
 \[\frac{
\mu\left(S(I)\right )\left(\log \frac{2\pi }{\vert I\vert }\right
)^\alpha }{|I|^s}\le C, \quad\hbox{for any interval
$I\subset\partial\D $}.
\]
\par
A positive Borel measure $\mu $ on $[0, 1)$ can be seen as a Borel
measure on $\mathbb D$ by identifying it with the measure $\tilde
\mu $ defined by $$ \tilde \mu (A)\,=\,\mu \left (A\cap [0,1)\right
),\quad \text{for any Borel subset $A$ of $\mathbb D$}.$$  In this
way a positive Borel measure $\mu $ on $[0, 1)$ is an $s$-Carleson
measure if and only if there exists a positive constant $C$ such
that
\[
\mu\left([t,1)\right )\le C(1-t)^s, \quad 0\le t<1.
\]
We have a similar statement for $\alpha$-logarithmic $s$-Carleson
measures.
\par The action of the operators $\,\mathcal I_\mu \,$ and $\,\mathcal H_\mu \,$
on distinct spaces of analytic functions have been studied in a
number of articles (see, e.\,\@g.,
\cite{Bao-Wu,Ch-Gi-Pe,Ga-Pe2010,GM1,GM2,Mer,Pow,Wi}).
\par Combining results of \cite{Ga-Pe2010} and of \cite{GM2} we can
state the following result.
\begin{other}\label{HmuH1} Let $\,\mu \,$ be a finite positive Borel
measure on $\,[0,1)$.
\begin{itemize} \item[(i)] The operator $\,\mathcal I_\mu \,$ is
well defined on $\,H^1\,$ if and only if $\,\mu \,$ is a Carleson
measure.
\item[(ii)] If $\,\mu \,$ is a Carleson
measure, then the operator $\,\mathcal H_\mu \,$ is also well
defined on $\,H^1\,$ and $\,\mathcal I_\mu f=\mathcal H_\mu f\,$ for
all $\,f\in H^1$.
\item[(iii)] The operator $\,\mathcal H_\mu \,$ is a bounded
operator from $\,H^1\,$ into itself if and only if $\,\mu \,$ is a
$1$-logarithmic $1$-Carleson measure.
\end{itemize}
\end{other}
\par Galanopoulos and Pel\'{a}ez
\cite[Theorem\,\@2.\,\@2]{Ga-Pe2010} proved the following.
\begin{other}\label{Carleson-cima} Let $\,\mu \,$ be a positive
Borel measure on $\,[0,1)$. If $\,\mu \,$ is a Carleson measure then
$\,\mathcal H_\mu (H^1)\subset \mathscr C$.\end{other}

This result is stronger than Theorem\,\@\ref{th-cima}(i). In view of
these results, the following question arises naturally.

\begin{Question}\label{q2} Suppose that $\,\mu \,$ is a $1$-logarithmic
$1$-Carleson measure on $\,[0,1)$. What can we say about the image
$\,\mathcal H_\mu (H^1)\,$ of $\,H^1\,$ under the action of the
operator $\,\mathcal H_\mu $?
\end{Question}

\par To answer Question\,\@\ref{q2}, let us start noticing that it is
known that, for $0<p\le 2$, the space of Dirichlet type $\,\Dp$ is
continuously included in $H^p\,$ (see \cite[Lemma\,\@1.\,\@4]{Vi}).
In particular, the space $\,\mathcal D^1_0\,$ is continuously
included in $\,H^1$. In fact, the estimates obtained by Vinogradov
in the proof of his lemma easily yield the inequality
$$ \Vert f\Vert _{H^1}\le 2\Vert f\Vert _{\mathcal D^1_0},\quad f\in
\mathcal D^1_0.$$ We shall prove that if $\,\mu \,$ is a
$1$-logarithmic $1$-Carleson measure on $\,[0,1)\,$ then $\,\mathcal
H_\mu (H^1)\,$ is contained in the space $\,\mathcal D^1_0$.
Actually, we have the following stronger result.
\begin{theorem}\label{H1D1} Let $\,\mu \,$ be a positive Borel
measure on $\,[0,1)$. Then the following conditions are equivalent.
\begin{itemize}
\item[(i)] $\mu \,$ is a $1$-logarithmic
$1$-Carleson measure.
\item[(ii)] $\mathcal H_\mu \,$ is a bounded operator from $\,H^1\,$
into itself.
\item[(iii)] $\mathcal H_\mu \,$ is a bounded operator from $\,H^1\,$
into $\,\mathcal D^1_0$.
\item[(iv)] $\mathcal H_\mu \,$ is a bounded operator from $\,\mathcal D^1_0\,$
into $\,\mathcal D^1_0$.
\end{itemize}
\end{theorem}
\par There is a gap between Theorem\,\@\ref{Carleson-cima} and
Theorem\,\@\ref{H1D1} and so it is natural to discuss the range of
$\,H^1\,$ under the action of $\,\mathcal H_\mu \,$ when $\,\mu \,$
is an $\,\alpha $-logarithmic $1$-Carleson measure with $\,0<\alpha
<1$. We shall prove the following result.
\begin{theorem}\label{alpha-log} Let $\,\mu \,$ be a positive Borel
measure on $\,[0,1)$. Suppose that $\,0<\alpha <1\,$ and that $\,\mu
\,$ is an $\,\alpha $-logarithmic $1$-Carleson measure. Then
$\,\mathcal H_\mu \,$ maps $\,H^1\,$ into the space $\,\mathcal
D^1(\log ^{\alpha -1})\,$ defined as follows:
$$\mathcal D^1(\log ^{\alpha -1})=\left \{ f\in \hol (\D ) : \int_{\D} \vert
f^\prime (z)\vert \left (\log \frac{2}{1-\vert z\vert }\right
)^{\alpha -1}\,dA(z)<\infty \right \}.$$
 \end{theorem}
 \par\medskip These results will be proved in
 Section\,\@\ref{proofs}. Since the space of Dirichlet type $\,\mathcal D^1_0\,$ has showed up in a natural
 way in our work, it seems natural to study the action of the operators $\,\mathcal H_\mu \,$ and $\,\mathcal I_\mu
 \,$ on the Bergman spaces $\,A^p_{\alpha }\,$ and the Dirichlet spaces $\,\Dpa\,$
 for general values of the parameters $\,p\,$ and $\,\alpha $. This will be done in
 Section\,\@\ref{Ber-Dir}.
 \par Throughout this paper the letter $C$ denotes a positive constant that may
change from one step to the next. Moreover, for two real-valued
functions $E_1, E_2$ we write $E_1\lesssim E_2$, or $E_1\gtrsim
E_2$, if there exists a positive constant $C$ independent of the
arguments such that $E_1\leq C E_2$, respectively $E_1\ge C E_2$. If
we have $E_1\lesssim E_2$ and $E_1\gtrsim E_2$ simultaneously then
we say that $E_1$ and $E_2$ are equivalent and we write $E_1\asymp
E_2$.
\section{Proofs of the theorems\,\@\ref{H1D1} and \,\@\ref{alpha-log}}\label{proofs}

\begin{Pf}{\,\em{Theorem\,\@\ref{H1D1}.}} We already know that (i)
and (ii) are equivalent by Theorem\,\@\ref{HmuH1}. \par To prove
that (i) implies (iii) we shall use some results about the Bloch
space. We recall that a function $\,f\in \hol (\mathbb D)\,$ is said
to be a Bloch function if
$$\Vert f\Vert _{\mathcal B}\,\ig \,\vert f(0)\vert \,+\,\sup_{z\in
\mathbb D}(1-\vert z\vert ^2)\vert f^\prime (z)\vert \,<\,\infty .$$
The space of all Bloch functions will be denoted by $\,\mathcal B$.
It is a non-separable Banach space with the norm $\,\Vert \cdot
\Vert _{\mathcal B}\,$ just defined. A classical source for the
theory of Bloch functions is \cite{ACP}. The closure of the
polynomials in the Bloch norm is the {\it {little Bloch space}}
$\,\mathcal B_0\,$ which consists of those $\,f\in \hol (\mathbb
D)\,$ with the property that
\begin{equation*}\lim_{\vert z\vert \to 1}(1-\vert z\vert ^2)\vert
f^\prime (z)\vert =0.\end{equation*} It is well known that (see
\cite[p.\,\@13]{ACP})
\begin{equation}\label{gr-bloch}\vert f(z)\vert \lesssim \Vert
f\Vert _{\mathcal B}\log \frac{2}{1-\vert z\vert }.\end{equation}
\par
The basic ingredient to prove that (i) implies (iii) is the fact
that the dual $\,\left (\mathcal B_0\right )^*\,$ of the little
Bloch space can be identified with the Bergman space $\,A^1$ via the
integral pairing
\begin{equation}\label{pairing-B0}
\langle h,f\rangle\,=\,\int _{\mathbb D}\,h(z)\,\overline
{f(z)}\,dA(z),\quad h\in \mathcal B_0, \, f\in A^1.
\end{equation}
(See \cite[Theorem\,\@5.\,\@15]{Zhu}).
\par Let us proceed to prove the implication (i)\,$\Rightarrow
$\,(iii). Assume that $\,\mu\,$ is a $1$-logarithmic $1$-Carleson
measure and take $\,f\in H^1$. We have to show that $\,\mathcal
I_\mu f\in \mathcal D^1_0\,$ or, equivalently, that $\,\left
(\mathcal I_\mu f\right )^\prime \in A^1$. Since $\,\mathcal B_0\,$
is the closure of the polynomials in the Bloch norm, it suffices to
show that
\begin{equation}\label{in-pol}\left \vert \int_{\mathbb D}\,h(z)\,\overline {\left (\mathcal I_\mu f\right
)^\prime (z)}\,dA(z)\right \vert \,\lesssim \Vert h\Vert_{\mathcal
B}\Vert f\Vert_{H^1},\quad\text{for any polynomial
$h$}.\end{equation}
\par
So, let $\,h\,$ be a polynomial. We have
\begin{align*} \int_{\mathbb D}\,h(z)\,\overline {\left (\mathcal I_\mu f\right
)^\prime (z)}\,dA(z)\,=\,&
 \int_{\mathbb D}\,h(z)\,\overline{\left
 (\int_{[0,1)}\frac{t\,f(t)}{(1-tz)^2}\,d\mu (t)\right )}\,dA(z)\\=\,&
\int_{\mathbb D}\,h(z)\,\int_{[0,1)}\frac{t\,\overline
{f(t)}}{(1-t\,\overline{z})^2}\,d\mu (t)\,dA(z)\\=\,&
\int_{[0,1)}\,t\,\overline {f(t)}\int_{\mathbb
D}\,\frac{h(z)}{(1-t\,\overline {z})^2}\,dA(z)\,d\mu (t).
\end{align*}
Because of the reproducing property of the Bergman kernel
\cite[Proposition\,\@4.\,\@23]{Zhu},  $\,\int_{\mathbb
D}\,\frac{h(z)}{(1-t\,\overline {z})^2}\,dA(z)\,=\,h(t)$. Then it
follows that \begin{equation}\label{ddd}\int_{\mathbb
D}\,h(z)\,\overline {\left (\mathcal I_\mu f\right )^\prime
(z)}\,dA(z)\,=\,\int_{[0,1)}\,t\,\overline {f(t)}\,h(t)\,d\mu
(t).\end{equation} Since $\,\mu \,$ is a $1$-logarithmic
$1$-Carleson measure, the measure $\,\nu \,$ defined by
$$d\nu (t)\,=\,\log \frac{2}{1-t}\,d\mu (t)$$ is a Carleson measure
\cite[Proposition\,\@2.\,\@5]{GM1}. This implies that
$$\int_{[0,1)}\vert f(t)\vert \log\frac{2}{1-t}\,d\mu (t)\lesssim
\Vert f\Vert _{H^1}.$$ This and (\ref{gr-bloch}) yield
$$\int_{[0,1)}\left\vert t\,\overline {f(t)}\,h(t)\right\vert\,d\mu (t)\,\lesssim \,\Vert
h\Vert_{\mathcal B}\Vert f\Vert _{H^1}.$$ Using this and
(\ref{ddd}), (\ref{in-pol}) follows.
\par
Since $\,\mathcal D^1_0\subset H^1\,$, the implication (iii)\,
$\Rightarrow$\, (iv) is trivial. To prove that (iv) implies (i) we
shall use the following result of Pavlovi\'c
\cite[Theorem\,\@3.\,\@2]{Pav-dec}.
\begin{other}\label{Pav-dec-D1} Let $\,f\in \hol (\D)$,
$\,f(z)=\sum_{n=0}^\infty a_nz^n$, and suppose that the sequence
$\,\{ a_n\} $ is a decreasing sequence of non-negative real numbers.
Then $\,f\in \mathcal D^1_0\,$ if and only if $\,\sum_{n=0}^\infty
\frac{a_n}{n+1}<\infty $, and we have
$$\Vert f\Vert _{\mathcal D^1_0}\,\asymp \,\sum_{n=0}^\infty
\frac{a_n}{n+1}.$$
\end{other}
\par
Now we turn to prove the implication (iv)\,$\Rightarrow $\,(i).
Assume that $\mathcal H_\mu $ is a bounded operator from $\,\mathcal
D^1_0\,$ into $\,\mathcal D^1_0$. We argue as in the proof of
Theorem\,\@1.\,\@1 of \cite{GM2}. For $\,\frac{1}{2}<b<1\,$ set
$$f_b(z)=\frac{1-b^2}{(1-bz)^2},\quad z\in \mathbb D.$$
We have $\,f_b^\prime (z)=\frac{2b(1-b^2)}{(1-bz)^3}\,$ ($z\in \D$).
Then, using Lemma\,\@3.\,\@10 of \cite{Zhu} with $\,t=0\,$ and
$\,c=1$, we see that
$$\Vert f_b\Vert _{\mathcal D^1_0}\,\asymp\,\int_{\D
}\frac{1-b^2}{\vert 1-bz\vert^3}\,dA(z)\,\asymp 1.$$ Since
$\,\mathcal H_\mu \,$ is bounded on $\,\mathcal D^1_0$, this implies
that
\begin{equation}\label{H1fb}1\gtrsim \Vert \mathcal H_\mu (f_b)\Vert
_{\mathcal D^1_0}.\end{equation} We also have,
$$f_b(z)\,=\,\sum_{k=0}^\infty a_{k,b}z^k,\quad\text{with
$a_{k,b}=(1-b^2)(k+1)b^k$.}$$ Since the $a_{k,b}$'s are positive, it
is clear that the sequence $\{ \sum_{k=0}^\infty \mu_{n+k}a_{k,b}\}
_{n=0}^\infty $ of the Taylor coefficients of $\mathcal H_\mu (f_b)$
is a decreasing sequence of non-negative real numbers. Using this,
Theorem\,\@\ref{Pav-dec-D1}, (\ref{H1fb}), and the definition of the
$a_{k,b}$'s, we obtain
\begin{align*} 1\,&\gtrsim \,\Vert \mathcal H_\mu (f_b)\Vert
_{\mathcal D^1_0}\,\gtrsim \,\sum_{n=1}^\infty \frac{1}{n}\left
(\sum_{k=0}^\infty \mu_{n+k}a_{k,b}\right )\\ & =  \sum_{n=1}^\infty
\frac{1}{n}\left (\sum_{k=0}^\infty
a_{k,b}\int_{[0,1)}t^{n+k}\,d\mu (t)\right )\\
& \gtrsim \,(1-b^2)\sum_{n=1}^\infty \frac{1}{n}\left
(\sum_{k=1}^\infty kb^k\int_{[b,1)}t^{n+k}\,d\mu (t)\right )\\
& \gtrsim \,(1-b^2)\sum_{n=1}^\infty \frac{1}{n}\left
(\sum_{k=1}^\infty kb^{n+2k}\,\mu\left ([b,1)\right )\right )
\\ &
=\,(1-b^2)\mu \left ([b,1)\right )\sum_{n=1}^\infty
\frac{b^n}{n}\left (\sum_{k=1}^\infty kb^{2k}\right )
\\ & =
\,(1-b^2)\mu \left ([b,1)\right )\left (\log \frac{1}{1-b}\right
)\frac{b^2}{(1-b^2)^2}.
\end{align*}
Then it follows that
$$\mu \left ([b,1)\right )\,=\,\og \left (
\frac{1-b}{\log\frac{1}{1-b}}\right ),\quad \text{as $b\to 1$}.$$
Hence, $\mu $ is a $1$-logarithmic $1$-Carleson measure.
\end{Pf}
\par\medskip
Before embarking on the proof of Theorem\,\@\ref{alpha-log} we have
to introduce some notation and results. Following \cite{Pav-dual},
for $\,\alpha \in \mathbb R\,$ the weighted Bergman space
$\,A^1(\log^\alpha )\,$
  consists of
those $\,f\in\hol (\D )\,$ such that
$$\Vert f\Vert _{A^1(\log^\alpha )}\,\ig \,\int_{\D }\vert f(z)\vert
\,\left (\log\frac{2}{1-\vert z\vert }\right )^\alpha\,dA(z)<\infty
.$$ This is a Banach space with the norm
$\,\Vert\cdot\Vert_{A^1(\log^\alpha )}\,$ just defined and the
polynomials are dense in $\,A^1(\log^\alpha )$. Likewise, we define
$$\mathcal D^1(\log^\alpha )=\{ f\in \hol (\D ) : f^\prime \in
A^1(\log^\alpha )\} .$$ \par We define also the Bloch-type space
$\,\mathcal B(\log^\alpha )\,$ as the space of those $\,f\in \hol
(\D )\,$ such that
$$\Vert f\Vert _{\mathcal B(\log^\alpha )}\ig \vert f(0)\vert
+\sup_{z\in \mathbb D}(1-\vert z\vert ^2)\left (\log
\frac{2}{1-\vert z\vert }\right )^{-\alpha }\vert f^\prime (z)\vert
<\infty ,$$ and  $$\mathcal B_0(\log^\alpha )=\left \{ f\in \hol (\D
) : \vert f^\prime (z)\vert =\op \left (\frac{\left
(\log\frac{2}{1-\vert z\vert }\right )^\alpha }{1-\vert z\vert
}\right ) ,\,\,\,\text{as $\,\vert z\vert \to 1$}\right \} .$$ The
space $\,\mathcal B(\log^\alpha )\,$ is a Banach space and
$\,\mathcal B_0(\log^\alpha )\,$ is the closure of the polynomials
in $\,\mathcal B(\log^\alpha )$.
\par We remark that the spaces $\,\mathcal D^1(\log^\alpha )\,$, $\,\mathcal B(\log^\alpha
)\,$, and $\,\mathcal B_0(\log^\alpha )\,$ were called $\,\mathfrak
B^1_{\log ^\alpha }$, $\,\mathfrak B_{\log ^\alpha }$, and
$\,\mathfrak b_{\log ^\alpha }\,$ in \cite{Pav-dual}. Pavlovi\'c
identified in \cite[Theorem\,\@2.\,\@4]{Pav-dual} the dual of the
space $\,\mathcal B_0(\log^\alpha )$.
\begin{other}\label{duality} Let $\,\alpha \in \mathbb R$. Then the
dual of $\,\mathcal B_0(\log^\alpha )\,$ is $\,A^1(\log^\alpha )\,$
via the pairing
$$\langle h, g\rangle\,=\,\int_{\mathbb D}\,f(z)\,\overline
{g(z)}\,dA(z),\quad h\in  \mathcal B_0(\log^\alpha ),\,\,\,g \in
A^1(\log^\alpha ).$$
\end{other}
\par Actually, Pavlovi\'c formulated the duality theorem in another
way but it is a simple exercise to show that his formulation is
equivalent to this one which is better suited to our work.
\par\medskip
\begin{Pf}{\,\em{Theorem\,\@\ref{alpha-log}.}} Let $\,\mu\,$ be a
positive Borel measure on $\,[0,1)\,$ and $\,0<\alpha <1$. Suppose
that $\,\mu \,$ is an $\,\alpha $-logarithmic $1$-Carleson measure.
Take $\,f\in H^1$. We have to show that $\,\mathcal I_\mu f\in
\mathcal D^1(\log ^{\alpha -1})\,$ or, equivalently, that $\,\left
(\mathcal I_\mu f\right )^\prime \in A^1(\log ^{\alpha -1})$.
Bearing in mind Theorem\,\@\ref{duality} and the fact that
$\,\mathcal B_0(\log ^{\alpha -1})\,$ is the closure of the
polynomials in $\,\mathcal B(\log ^{\alpha -1})\,$, it suffices to
show that
\begin{equation}\label{in-pol-2}\left \vert \int_{\mathbb D}\,h(z)\,\overline {\left (\mathcal I_\mu f\right
)^\prime (z)}\,dA(z)\right \vert \,\lesssim \Vert h\Vert_{\mathcal
B(\log^{\alpha -1})}\Vert f\Vert_{H^1},\quad\text{for any polynomial
$h$}.\end{equation}
\par
So, let $\,h\,$ be a polynomial. Arguing as in the proof of the
implication (i) $\Rightarrow $ (iii) in Theorem\,\@\ref{H1D1} we
obtain \begin{equation}\label{pre}\int_{\D }h(z)\,\overline {\left
(\mathcal I_\mu f\right )^\prime
(z)}\,dA(z)\,=\,\int_{[0,1)}t\,\overline {f(t)}\,h(t)\,d\mu
(t).\end{equation} Now, it is clear that
$$\vert h(z)\vert \lesssim \Vert h\Vert _{\mathcal B(\log ^{\alpha
-1})}\left (\log \frac{2}{1-\vert z\vert }\right )^\alpha ,$$ and
then it follows that
$$\int_{[0,1)}\left \vert t\,\overline {f(t)}\,h(t)\right\vert \,d\mu
(t)\lesssim \Vert h\Vert _{\mathcal B(\log ^{\alpha
-1})}\int_{[0,1)}\vert f(t)\vert \left (\log\frac{2}{1-t}\right
)^{\alpha }\,d\mu (t).$$ Using the fact that the measure $\,\left
(\log\frac{2}{1-t}\right )^{\alpha }\,d\mu (t)\,$ is a Carleson
measure \cite[Proposition\,\@2.\,\@5]{GM1}, this implies that
$$\int_{[0,1)}\left \vert t\,\overline {f(t)}\,h(t)\right\vert \,d\mu
(t)\lesssim \Vert h\Vert _{\mathcal B(\log ^{\alpha -1})}\Vert
f\Vert _{H^1}.$$ This and (\ref{pre}) give (\ref{in-pol-2}).
\end{Pf}

\section{The operators $\mathcal H_\mu $ acting on  Bergman spaces and on Dirichlet spaces}\label{Ber-Dir}
Jevti\'c and Karapetrovi\'c \cite{JK} have recently proved the
following result.
\begin{other}\label{JK-H-Dir} The Hilbert operator $\,\mathcal H\,$
is a bounded operator from $\,\Dpa \,$ into itself if and only if
$\,\max (-1, p-2)\,<\alpha \,<2p-2$.\end{other}
\par Now, it is well known that $\,A^p_\alpha \,=\,\mathcal
D^p_{\alpha +p}\,$ (see \cite[Theorem\,\@4.\,\@28]{Zhu}). Hence,
regarding Bergman spaces Theorem\,\@\ref{JK-H-Dir} says the
following. \begin{otherc}\label{H-bounded-Ber} The Hilbert operator
$\,\mathcal H\,$ is a bounded operator from $\,A^p_\alpha\,$ into
itself if and only if $\,-1<\alpha <p-2$.\end{otherc}
\par Let us recall that Diamantopoulos \cite{Di} had proved before that the Hilbert operator is
bounded on $\,A^p\,$ for $\,p>2$, but not on $\,A^2$. The situation
on $\,A^2\,$ is even worse. Dostani\'c, Jevti\'c, and Vukoti\'c
\cite{DJV} proved that the Hilbert operator is not well defined on
$\,A^2$. Indeed, they considered the function $\,f\,$ defined by
\begin{equation}\label{FApp-2}f(z)\,=\,\sum_{n=1}^\infty \frac{1}{\log (n+1)}\,z^n,\quad z\in
\D ,\end{equation} which belongs to $\,A^2.$ However, the series
defining $\,\mathcal Hf(0)\,$ is $\,\sum_{n=1}^\infty
\frac{1}{(n+1)\log (n+1)}=\infty \,$ and the integral defining
$\,\mathcal If(0)\,$ is $\,\int_0^1\,f(t)\,dt\,=\,\infty $. Hence
neither $\,\mathcal H\,$ nor $\,\mathcal I\,$ are defined on
$\,A^2$.
\par This result can be extended.
We can assert that $\,\mathcal H\,$ is not well defined on
$\,A^p_{p-2}\,$ for any $\,p>1$. Indeed, let $\,f\,$ be the function
defined in (\ref{FApp-2}). Notice that the sequence $\,\{
\frac{1}{(n+1)\log (n+1)}\} $ is decreasing and that
$\,\sum_{n=1}^\infty \frac{1}{n\left (\log (n+1)\right )^p}<\infty
$. Then (see Proposition\,\@\ref{Apa-dec} below) it follows that
$\,f\in A^p_{p-2}$, and we have already seen that $\,\mathcal Hf\,$
and $\,\mathcal If\,$ are not defined. Since $\,\alpha \ge
p-2\,\,\Rightarrow\,\, A^p_{p-2} \subset A^p_\alpha $, it follows
that the Hilbert operator $\,\mathcal H\,$ is not defined on
$\,A^p_\alpha \,$ if $\,\alpha \ge p-2$.
\par In this section we shall obtain extensions of the mentioned results of Jevti\'c and Karapetrovi\'c
considering the generalized Hilbert operators $\,\mathcal H_\mu \,$.
\begin{theorem}\label{Hmu-carl-well-Dir}
Suppose that $\,\max (-1, p-2)\,<\alpha \,<2p-2\,$ and let $\,\mu
\,$ be a finite positive Borel measure on $\,[0,1)$. If $\,\mu \,$
is a Carleson measure then the operators $\,\mathcal H_\mu \,$ and
$\,\mathcal I_\mu \,$ are well defined on $\,\Dpa $. Furthermore,
$\mathcal I_\mu f=\mathcal H_\mu f$, for all $\,f\in \Dpa $.
\end{theorem}
\par When dealing with Bergman spaces Theorem\,\@\ref{Hmu-carl-well-Dir}
reduces to the following.
\begin{corollary}\label{Hmu-carl-well-Ber}Suppose that $\,p>1\,$ and $\,-1\,<\alpha\,<\,p-2\,$, and let $\,\mu
\,$ be a finite positive Borel measure on $\,[0,1)$. If $\,\mu \,$
is a Carleson measure then the operators $\,\mathcal H_\mu \,$ and
$\,\mathcal I_\mu \,$ are well defined on $\,A^p_\alpha $.
Furthermore, $\mathcal I_\mu f=\mathcal H_\mu f$, for all $\,f\in
A^p_\alpha $.\end{corollary}
\begin{Pf}{\,\em{Theorem\,\@\ref{Hmu-carl-well-Dir}.}}
Suppose that $\,\mu \,$ is a Carleson measure and take $\,f\in \Dpa
$. Set $\,\beta =\frac{2+\alpha }{p}-1$. Observe that $\,0<\beta
<1$. Using \cite[Theorem\,\@4.\,\@14]{Zhu}, we see that $\,\vert
f^\prime (z)\vert \,\lesssim \frac{1}{(1-\vert z\vert )^{(2+\alpha
)/p}}\,$ and, hence, $\,\vert f(z)\vert \,\lesssim \frac{1}{(1-\vert
z\vert )^\beta }$. Then it follows that
\begin{equation*}\int_{[0,1)}\vert f(t)\vert \,d\mu
(t)\,\lesssim \,\int_{[0,1)}\frac{d\mu (t)}{(1-t)^\beta
}.\end{equation*} Integrating by parts, using that $\,\mu \,$ is a
Carleson measure, and that $\,0<\beta <1$, we obtain
\begin{align*}
\int_{[0,1)}\frac{d\mu (t)}{(1-t)^\beta }\,&=\,\mu
([0,1))\,-\,\lim_{t\to 1}\frac{\mu ([t,1))}{(1-t)^{\beta
}}\,+\,\beta \int_0^1\,\frac{\mu ([t,1))}{(1-t)^{\beta +1}}\,dt
\\ &=\,\mu
([0,1))\,+\,\beta \int_0^1\,\frac{\mu ([t,1))}{(1-t)^{\beta +1}}\,dt
\\ &\lesssim \,\mu
([0,1))\,+\,\int_0^1\,\frac{1}{(1-t)^{\beta }}\,dt \\ &<\,\infty .
\end{align*}
Consequently, we obtain that
\begin{equation}\label{int-fin}\int_{[0,1)}\vert f(t)\vert \,d\mu
(t)\,<\,\infty .\end{equation} Clearly, this implies that the
integral
\begin{equation}\label{uni-abs}\int_{[0,1)}\frac{f(t)\,d\mu
(t)}{1-tz}\quad\text{converges absolutely and uniformly on compact
subsets of $\D$}.\end{equation} This gives that $\,\mathcal I_\mu
f\,$ is a well defined analytic function in $\D $ and that
\begin{equation}\label{Isum}\mathcal I_\mu f(z)=\sum _{n=0}^\infty
\left (\int_{[0,1)}t^n\,f(t)\,d\mu (t)\right )z^n,\quad z\in
\D.\end{equation}
\par Using \cite[Theorem\,\@2.\,\@1]{JK-Fil} (see also
\cite[Theorem\,\@2.\,\@1]{JK}) we see that for these values of
$\,p\,$ and $\,\alpha \,$ we have that if $\,f\in A^p_\alpha $,
$f(z)=\sum_{n=0}^\infty a_nz^n$, then $\,\sum_{k=0}^\infty
\frac{\vert a_k\vert }{k+1}<\infty $. Now, since $\,\mu \,$ is a
Carleson measure we have that $\,\vert \mu _n\vert \lesssim
\frac{1}{n+1}$ (\cite[Proposition\,\@1]{Ch-Gi-Pe}). Then it follows
that
$$\sum_{k=0}^\infty \vert \mu_{n+k}a_k\vert \,\lesssim \,\sum_{k=0}^\infty
\frac{\vert a_k\vert }{k+n+1}\lesssim  \sum_{k=0}^\infty \frac{\vert
a_k\vert }{k+1},\quad\text{for all $\,n$}.$$ Clearly, this implies
that $\,\mathcal H_\mu f\,$ is a well defined analytic function in
$\,\D\,$ and that $\,\int_{[0,1)}t^n\,f(t)\,d\mu
(t)\,=\,\sum_{k=0}^\infty \mu_{n+k}a_k\,$ for all $\,n$. This and
(\ref{uni-abs}) give that $\,\mathcal I_\mu f\,=\,\mathcal H_\mu f$.
\end{Pf}
\par\medskip
Our next result is an extension of Corollary\,\@\ref{H-bounded-Ber}
\begin{theorem}\label{Hmu-bounded-Ber} Suppose that $\,-1\,<\alpha\,<\,p-2\,$ and let $\,\mu
\,$ be a finite positive Borel measure on $\,[0,1)$. \par The
operator $\,\mathcal H_\mu \,$ is well defined on $\,A^p_\alpha \,$
and it is a bounded operator from $\,A^p_\alpha \,$ to itself if and
only if $\,\mu \,$ is a Carleson measure.
\end{theorem}
\par\medskip A number of results will be needed to prove
this theorem. We start with a characterization of the functions
$\,f\in \hol (\D )\,$ whose sequence of Taylor coefficients is
decreasing which belong to $\,A^p_\alpha $.
\begin{proposition}\label{Apa-dec} Let $\,f\in \hol (\D )$, $f(z)=\sum_{n=0}^\infty
a_n\,z^n$ $\,(z\in \D)$. Suppose that $\,1<p<\infty $, $\,\alpha
>-1$, and that the sequence $\{ a_n\} _{n=0}^\infty $ is a
decreasing sequence of non-negative real numbers. Then
$$f\in A^p_\alpha \,\,\,\Leftrightarrow\,\,\,\sum_{n=1}^\infty
n^{p-3-\alpha }a_n^p\,<\,\infty .$$ Furthermore, $\Vert
f\Vert_{A^p_\alpha }^p\,\asymp\,\vert
a_0\vert^p\,+\,\sum_{n=1}^\infty n^{p-3-\alpha }a_n^p\,<\,\infty $.
\end{proposition}
\par This result can be proved with arguments similar to those used
in the proofs of \cite[Theorem\,3.\@10]{GM1} and
\cite[Theorem\,\@3.\,\@1]{Pav-dec} where the analogous results for
the Besov spaces $\,B^p=\mathcal D^p_{p-2}$ ($p>1$) and for the
spaces $\,\Dp\,$ ($p>1$) were proved. The case $\,\alpha =0\,$ is
proved in \cite[Proposition\,\@2.\,\@4]{BKV}. Consequently, we omit
the details.
\par The following lemma is a generalization of
\cite[Lemma\,\@3\,\@(ii)]{GaGiPeSis}.
\begin{lemma}\label{le-infty-Apa} Let $\,\mu \,$ be a positive Borel measure on $\,[0,1)\,$ which is a Carleson measure.
Assume that $\,0<p<\infty\,$ and $\,\alpha >-1$. Then there exists a
positive constant $\,C=C(p,\alpha ,\mu )\,$ such that for any
$\,f\in A^p_\alpha \,$
$$\int_{[0,1)}M_\infty ^p(r,f)(1-r)^{\alpha +1}\,d\mu (r)\,\le C\Vert
f\Vert _{A^p_\alpha }^p.$$
\end{lemma}
\par Of course, $M_\infty (r,f)=\sup_{\vert z\vert =r}\vert
f(z)\vert $.\par \begin{pf} Take $\,f\in A^p_\alpha \,$ and set
$$g(r)=M_\infty^p(r,f)(1-r)^{\alpha +1},\,\, F(r)=\mu ([0,r))\,-\mu ([0,1))\,=\,-\mu ([r,1)),\,\, 0<r<1.$$
Integrating by parts, we have
\begin{align}\label{aaaa}& \int_{[0,1)}M_\infty ^p(r,f)(1-r)^{\alpha +1}\,d\mu
(r)\,=\,\int_{[0,1)}g(r)\,d\mu (r)\\ \,=&\,\lim_{r\to
1}g(r)F(r)\,-\,g(0)F(0)\,-\int_0^1g^\prime (r)F(r)\,dr \nonumber\\
\,=&\,\vert f(0)\vert ^p\mu ([0,1))\,-\,\lim_{r\to 1}M_\infty
^p(r,f)(1-r)^{\alpha +1}\mu([r,1))\,+\,\int_0^1g^\prime (r)\mu ([r,
1))\,dr.\nonumber
\end{align} Since $\,f\in A^p_\alpha \,$ we have that $\,M_\infty
^p(r,f)=\op \left ((1-r)^{-2-\alpha }\right )$, as $\,r\to 1$ (see,
e.\,\@g., \cite[p.\,\@54]{HKZ}). This and the fact that $\,\mu \,$
is a Carleson measure imply that
\begin{equation}\label{limit-zero}\lim_{r\to 1}M_\infty
^p(r,f)(1-r)^{\alpha +1}\mu([r,1))\,=\,0.\end{equation} Using again
that $\,\mu \,$ is a Carleson measure and integrating by parts we
see that
\begin{align*}\int_0^1g^\prime
(r)\mu ([r, 1))\,dr\,\lesssim \,&\int_0^1g^\prime (r)(1-r)\,dr\\
\,=\,&\lim_{r\to 1} g(r)(1-r)\,-g(0)\,+\,\int_0^1g(r)\,dr\\
\,\le \,& \lim_{r\to 1}M_\infty ^p(r,f)(1-r)^{\alpha
+2}\,+\,\int_0^1M_\infty ^p(r,f)(1-r)^{\alpha +1}\,dr \\
\,=\,&\int_0^1M_\infty ^p(r,f)(1-r)^{\alpha +1}\,dr.
\end{align*}
Then, using \cite[Lemma \,\@3.\,\@(ii)]{GaGiPeSis}, it follows that
\begin{align*}\int_0^1g^\prime
(r)\mu ([r, 1))\,dr\,\lesssim \,\Vert f\Vert_{A^p_\alpha
}^p.\end{align*} Using this and (\ref{limit-zero}) in (\ref{aaaa})
readily yields $\,\int_{[0,1)}M_\infty ^p(r,f)(1-r)^{\alpha
+1}\,d\mu (r)\,\lesssim \,\Vert f\Vert_{A^p_\alpha }^p$.
\end{pf}
\par
We shall also need the following characterization of the dual of the
spaces $\,A^q_\beta $ ($q>1$). It is a special case of
\cite[Theorem\,\@2.\,\@1]{Lu}.
\begin{lemma}\label{dualApa}  If  $1<q<\infty$ and $\beta >-1$, then  the dual of $A^q_\beta $
can be identified with  $A^{p}_{\alpha}$ where
$\frac{1}{p}+\frac{1}{q}=1$ and $\alpha $ is any number with $\alpha
>-1$, under the pairing
\begin{equation}\label{pair-Berg}
\langle h, f\rangle_{A_{q,\beta ,\alpha }}=\int_\D
h(z)\overline{f(z)}(1-|z|^2)^{\frac{\beta }{q}+\frac{\alpha
}{p}}\,dA(z),\quad h\in A^q_\beta ,\,\,\,\,f\in A^p_\alpha .
\end{equation}
\end{lemma}
\par
Finally, we recall the following result from \cite[(5.\,\@2),
p.\,\@242]{GaGiPeSis} which is a version of the classical Hardy's
inequality \cite[pp.\,\@244-245]{HLP}.
\begin{lemma}\label{Hardy-ineq} Suppose that $\,k>0$, $\,q>1$, and
$\,h\,$ is a non-negative function defined in $\,(0,1)$, then
\begin{equation*} \int_0^1\left(\int_{1-r}^1 h(t)\, dt \right)^q (1-r)^{k-1}\,dr \leq
\left(\frac{q}{k}\right)^q \int_0^1  (h(1-r))^q (1-r)^{q+k-1}\,dr.
 \end{equation*}
\end{lemma}
\par\medskip
\begin{Pf}{\,\em{Theorem\,\@\ref{Hmu-bounded-Ber}.}}
Suppose first that $\,\mathcal H_\mu\,$ is a bounded operator from
$\,A^p_\alpha \,$ into itself. For $0<b<1$, set
$$f_b(z)=\frac{(1-b^2)^{1-\frac{\alpha}{p}}}{(1-bz)^{\frac{2}{p}+1}},\quad z\in\D.$$
Recall that $\,p-\alpha >2$. Then using
\cite[Lemma\,\@3.\,\@10]{Zhu} with $\,t=\alpha \,$ and $\,c=p-\alpha
$, we obtain
\begin{align*}
\Vert f_b\Vert_{A^p_\alpha }^p\,=\,(1-b^2)^{p-\alpha }\int_{\D
}\frac{(1-\vert z\vert ^2)^\alpha }{\vert 1-bz\vert
^{2+p}}\,dA(z)\,\asymp \,1.\end{align*} Since $\,\mathcal H_\mu\,$
is bounded on $\,A^p_\alpha $, this implies
\begin{align}\label{fb1}1\,\gtrsim \Vert\mathcal H\mu (f_b)\Vert_{A^p_\alpha
}.\end{align} We also have $$f_b(z)=\sum_{k=0}^\infty
a_{k,b}z^k,\,\,\,(z\in \D ),\quad\text{with $\,a_{k,b}\asymp
(1-b^2)^{1-\frac{\alpha }{p}}k^{\frac{2}{p}}b^k.$}$$ Since the
$a_{k,b}$'s are positive, it is clear that the sequence $\left \{
\sum_{k=0}^\infty \mu_{n+k}a_{k,b}\right \} _{n=0}^\infty $ of the
Taylor coefficients of $\mathcal H_\mu (f_b)$ is a decreasing
sequence of non-negative real numbers. Using this,
Proposition\,\@\ref{Apa-dec}, (\ref{fb1}),
 and the definition of the $a_{k,b}$'s,
we obtain
\begin{align*}
1&\gtrsim \|\mathcal H_\mu(f_b)\|_{A^p_\alpha }^p\gtrsim
\sum_{n=1}^\infty n^{p-\alpha-3}\left(\sum_{k=1}^\infty
\mu_{n+k}a_{k,b}\right)^p
\\ & = \sum_{n=1}^\infty n^{p-\alpha-3}\left(\sum_{k=1}^\infty a_{k,b}\int_{[0,1)} t^{n+k} d\mu(t)\right)^p
\\ & \gtrsim (1-b^2)^{p-\alpha}\sum_{n=1}^\infty n^{p-\alpha-3}\left(\sum_{k=1}^\infty k^{\frac{2}{p}}b^{k}\int_{[b,1)} t^{n+k} d\mu(t)\right)^p
\\ & \ge (1-b^2)^{p-\alpha}\sum_{n=1}^\infty n^{p-\alpha-3}\left(\sum_{k=1}^\infty k^{\frac{2}{p}}b^{n+2k}\mu([b,1))\right)^p
\\ & = (1-b^2)^{p-\alpha}\mu([b,1))^p\sum_{n=1}^\infty n^{p-\alpha-3}b^{np}\left(\sum_{k=1}^\infty k^{\frac{2}{p}}b^{2k}\right)^p
\\ &\asymp (1-b^2)^{p-\alpha}\mu([b,1))^p\frac{1}{(1-b^2)^{2+p}}\sum_{n=1}^\infty n^{p-\alpha-3}b^{np}
\\ &\asymp (1-b^2)^{p-\alpha}\mu([b,1))^p\frac{1}{(1-b^2)^{2+p}}\cdot\frac{1}{(1-b^2)^{p-\alpha -2}}
\\ &\asymp \mu([b,1))^p\frac{1}{(1-b)^{p}}.
\end{align*}
Then it follows that $$\mu\left ([b,1)\right )\,=\,\og \left
(1-b\right ),\quad\text{as $b\to 1$},$$ and, hence, $\mu $ is a
Carleson measure.\par We turn to prove the other implication. So,
suppose that $\,\mu \,$ is a Carleson measure and take $\,f\in
A^p_\alpha $. Let $\,q\,$ be defined by the relation
$\,\frac{1}{p}+\frac{1}{q}=1\,$ and take $\,\beta =\frac{-\alpha
q}{p}=\frac{-\alpha }{p-1}$. Observe that $\,\beta >-1\,$ and that
with this election of $\,\beta \,$ the weight in the pairing
(\ref{pair-Berg}) is identically equal to $1$. We have to show that
$\,\mathcal H_\mu f\in A^p_\alpha \,$ which is equal to $\,\left
(A^q_\beta \right )^*\,$ under the pairing $\,\langle \cdot ,\cdot
\rangle_{q,\beta ,\alpha }$. So take $\,h\in A^q_\beta $.
\begin{align*}\langle h,\mathcal H_\mu f
\rangle_{q,\beta ,\alpha }\,=&\,\int_{\mathbb D}h(z)\,\overline
{\mathcal H_\mu f(z)}\,dA(z)\\ =&\, \int_{[0,1)}\overline{f(t)}\left
(\int_{\D }\frac{h(z)}{1-t\,\overline{z}}\,dA(z)\right )\,d\mu (t)\\
=&\, \int_{[0,1)}\overline{f(t)}\left (\int_0^1\,\frac{r}{\pi
}\,\int_0^{2\pi } \frac{h(re^{i\theta })}{1-tre^{-i\theta
}}\,d\theta\,dr\right )\,d\mu (t)
\\
=&\, \int_{[0,1)}\overline{f(t)}\left (\int_0^1\,\left (\frac{r}{\pi
i}\,\int_{\vert \xi\vert =1} \frac{h(r\xi )}{\xi -tr}\,d\xi\right
)\,dr\right )\,d\mu (t)
\\
=&\, 2\int_{[0,1)}\overline{f(t)}\left (\int_0^1\,rh(r^2t)\,dr\right
)\,d\mu (t).
\end{align*}
Thus,
\begin{align*}\left\vert \langle h,\mathcal H_\mu f
\rangle_{q,\beta ,\alpha }\right\vert \le 2\int_0^1\vert f(t)\vert
G(t)\,d\mu (t),\end{align*} where $G(t)=\int_0^1 r\vert h(r^2t)\vert
\,dr$. Using H\"older's inequality we obtain,
 \begin{align*} &\int_{[0,1)}
f(t)G(t)\,d\mu(t)=\int_{[0,1)} |f(t)|(1-t)^{\frac{\alpha+1}{p}}
G(t)(1-t)^{-\frac{\alpha+1}{p}}\,d\mu(t)
\\ &\le \left(\int_{[0,1)} |f(t)|^p(1-t)^{\alpha+1}\,d\mu(t) \right)^{1/p}
\cdot\left( \int_{[0,1)}
G(t)^q(1-t)^{-\frac{q(\alpha+1)}{p}}\,d\mu(t)\right)^{1/q}.
\end{align*}
Lemma\,\@\ref{le-infty-Apa} implies that
$$\left(\int_{[0,1)} |f(t)|^p(1-t)^{\alpha+1}\,d\mu(t)
\right)^{1/p}\,\lesssim \,\Vert f\Vert_{A^p_\alpha }.$$ Next we will
show that \begin{align}\label{claim}\int_{[0,1)}
G(t)^q(1-t)^{-\frac{q(\alpha+1)}{p}}\,d\mu(t)\lesssim
\|h\|^q_{A^q_\beta }.\end{align} This will give that
$$\left\vert \langle h,\mathcal H_\mu f
\rangle_{q,\beta ,\alpha }\right \vert\,\lesssim \,\Vert
f\Vert_{A^p_\alpha }\cdot \|h\|^q_{A^q_\beta }.$$ By the duality
theorem, this implies that $\,\mathcal H_\mu f\in A^p_\alpha $.
\par Let us prove (\ref{claim}).
Observe first that if $0<t<1/2$ then $|h(r^2t)|\le
M_\infty(\frac{1}{2},h)$ for each $r\in(0,1)$, thus
$$ G(t)=\int_0^1 |h(r^2t)|r\,dr\,\le M_\infty\left(\frac{1}{2},h\right), \quad 0<t<1/2.$$
Clearly, this implies
\begin{align}\label{0unmedio}
\int_{[0,1/2)}
G(t)^q(1-t)^{-\frac{q(\alpha+1)}{p}}\,d\mu(t)\,\lesssim M_\infty
^q\left (\frac{1}{2},h\right )\,\lesssim \,\|h\|^q_{A^q_\beta }.
\end{align}
Notice that $\,-\frac{q(\alpha +1)}{p}\,=\,\frac{p-2-\alpha
}{p-1}-1\,>-1$. Making the change of variables $\,r^2t=s$, we obtain
$\int_0^1 r|h(r^2t)|\,dr=\frac{1}{2t}\int_0^t |h(s)|\,ds$ and,
hence,
\begin{align}\label{G1medio1}&
\int_{[1/2,1)}  G(t)^q(1-t)^{-\frac{q(\alpha+1)}{p}}\,d\mu(t)\\ &=
\int_{[1/2,1)}  \left( \int_0^1
|h(r^2t)|r\,dr\right)^q(1-t)^{-\frac{q(\alpha+1)}{p}}\,d\mu(t)\nonumber
\\ &=\int_{[1/2,1)} \frac{1}{(2t)^q} \left( \int_0^t |h(s)|\,ds\right)^q(1-t)^{-\frac{q(\alpha+1)}{p}}\,d\mu(t)\nonumber
\\ &\le \int_{[1/2,1)} \left( \int_0^t M_\infty(s,h)\,ds\right)^q(1-t)^{-\frac{q(\alpha+1)}{p}}\,d\mu(t)\nonumber
\\ &\le \int_{[0,1)} \left( \int_{1-t}^1 M_\infty(1-s,h)\,ds\right)^q(1-t)^{-\frac{q(\alpha+1)}{p}}\,d\mu(t)\nonumber
\end{align}
Let us call $H(t)= \left( \int_{1-t}^1
M_\infty(1-s,h)\,ds\right)^q(1-t)^{-\frac{q(\alpha+1)}{p}}$ for
$0\le t <1$. Integrating by parts we obtain the following
\begin{align}\label{Htdmut}
\int_{[0,1)} H(t)\,d\mu(t)=H(0)\mu([0,1))-\lim_{t\to 1^-}
H(t)\mu([t,1))+ \int_0^1 \mu([t,1))H^\prime(t)\,dt.
\end{align}
The first term is equal to $0$. Using the fact that $\mu$ is a
Carleson measure we have that
\begin{align*}
H(t)\mu([t,1)) &\lesssim (1-t)H(t)
\\ &= \left( \int_{1-t}^1 M_\infty(1-s,h)\,ds\right)^q(1-t)^{1-\frac{q(\alpha+1)}{p}}
\\ &=\left( \int_{0}^t M_\infty(s,h)\,ds\right)^q(1-t)^{1-\frac{q(\alpha+1)}{p}}
.
\end{align*}
Since $\,h\in A^q_\beta \,$ we have $\,M_\infty(t,h)=o\left(
(1-t)^{-\frac{\beta+2}{q}}\right )$, as $\,t\to 1$. Then, bearing in
mind that $\,\frac{\beta +2}{q}>1$, it follows that
\begin{align}\label{secondop}H(t)\mu([t,1))\,=\,\op \left ((1-t)^{-\beta -2+q}\cdot
(1-t)^{1-\frac{q(\alpha+1)}{p}}\right )=\op (1),\quad \text{as
$\,t\to 1$}.\end{align} Actually, we have also proved that
\begin{align}\label{sssecondop}(1-t)H(t)\,=\,\op (1),\quad \text{as $\,t\to
1$}.\end{align}
\par Using that $\mu$ is a Carleson measure, integrating by parts, and using the definition of $\,H\,$ and (\ref{sssecondop}), we obtain
\begin{align}\label{mutHprime}
\int_0^1 \mu([t,1))H^\prime(t)\,dt&\lesssim \int_0^1
(1-t)H^\prime(t)\,dt
\\&= \lim_{t\to 1}(1-t)H(t)\,-\,H(0)\,+\,\int_0^1 H(t)\,dt \nonumber
\\ &= \int_0^1 \left( \int_{1-t}^1
M_\infty(1-s,h)\,ds\right)^q(1-t)^{-\frac{q(\alpha+1)}{p}}\,dt.\nonumber
\end{align}
Now, using Lemma\,\@\ref{Hardy-ineq} and
\cite[Lemma\,\@3]{GaGiPeSis}, we see that
\begin{align*}
\int_0^1 \left( \int_{1-t}^1
M_\infty(1-s,h)\,ds\right)^q(1-t)^{-\frac{q(\alpha+1)}{p}}\,dt
\lesssim \int_0^1 M^q_\infty(t,h)(1-t)^{\alpha+1}\,dt \lesssim
\|h\|^q_{A^q_\beta }.
\end{align*}
Using this, (\ref{mutHprime}), (\ref{secondop}), (\ref{Htdmut}), and
(\ref{G1medio1}), it follows that
$$\int_{[1/2,1)}  G(t)^q(1-t)^{-\frac{q(\alpha+1)}{p}}\,d\mu(t)\,\lesssim
\|h\|^q_{A^q_\beta }.$$ This and (\ref{0unmedio}) yield
(\ref{claim}).
\end{Pf}
\par\medskip
Our final aim in this article is to find the analogue of
Theorem\,\@\ref{Hmu-bounded-Ber} for Dirichlet spaces. In other
words, we wish give an answer to the following question.
\begin{Question}\label{q3} If  $\,\max (-1, p-2)<\alpha <2p-2$, is it true that
$\,\mathcal H_\mu \,$ is a bounded operator from $\,\Dpa\,$ into
itself if and only if $\,\mu \,$ is a Carleson measure?
\end{Question} Since
$\,p-1<\alpha <2p-2\,$ implies that $\,\Dpa\,=\,A^p_{\alpha
-p},$\,\,\, Theorem\,\@\ref{Hmu-bounded-Ber} answers the question
affirmatively for these values of $\,p\,$ and $\,\alpha $. It
remains to consider the case $\,\max (-1, p-2)<\alpha \le p-1$. We
shall prove the following result which gives a positive answer to
Question\,\@\ref{q3} in the case $\,p>1$.
\begin{theorem}\label{Hmu-bounded-Dir} Suppose that $\,p>1\,$ and
$\,p-2<\alpha \le p-1$, and let $\,\mu \,$ be a finite positive
Borel measure on $\,[0,1)$. \par The operator $\,\mathcal H_\mu \,$
is well defined on $\,\Dpa\,$ and it is a bounded operator from
$\,\Dpa\,$ into itself if and only if $\,\mu \,$ is a Carleson
measure.
\end{theorem}
\par The following two lemmas will be needed in the proof of
Theorem\,\@\ref{Hmu-bounded-Dir}. The first one follows trivially
from Proposition\,\@\ref{Apa-dec}.
\begin{lemma}\label{lemDpa} Let $\,f\in
\hol (\D )$, $f(z)=\sum_{n=0}^\infty a_nz^n\,$ $\,(z\in \D)$.
Suppose that $1\,<p<\infty \,$ and $\,p-2<\alpha\le p-1$, and that
the sequence $\,\{ a_n\} _{n=0}^\infty \,$ is a decreasing sequence
of non-negative real numbers. Then
$$f\in\Dpa\,\,\,\Leftrightarrow\,\,\,\sum_{n=0}^\infty (n+1)^{2p-\alpha-3}a_n^p<\infty.$$
\end{lemma}
\par The following lemma is a generalization of
\cite[Lemma\,\@4]{GaGiPeSis}.
\begin{lemma}\label{le-infty-Dpa} Let $\,\mu \,$ be a positive Borel measure on $\,[0,1)\,$ which is a Carleson measure.
Assume that $\,0<p<\infty\,$ and $\,\alpha >-1$. Then there exists a
positive constant $\,C=C(p,\alpha ,\mu )\,$ such that for any
$\,f\in \Dpa\,$
$$\int_{[0,1)}M_\infty ^p(r,f)(1-r)^{\alpha -p+1}\,d\mu (r)\,\le C\Vert
f\Vert _{\Dpa}^p.$$
\end{lemma}
\begin{pf} We argue as in the proof of Lemma\,\@\ref{le-infty-Apa}. Take $\,f\in \Dpa\,$ and set
$$g(r)=M_\infty^p(r,f)(1-r)^{\alpha -p+1},\,\, F(r)=\mu ([0,r))\,-\mu ([0,1))\,=\,-\mu ([r,1)),\,\, 0<r<1.$$
Integrating by parts, we have
\begin{align}\label{aaaaD}& \int_{[0,1)}M_\infty ^p(r,f)(1-r)^{\alpha -p+1}\,d\mu
(r)\,=\,\int_{[0,1)}g(r)\,d\mu (r)\\ \,=&\,\lim_{r\to
1}g(r)F(r)\,-\,g(0)F(0)\,-\int_0^1g^\prime (r)F(r)\,dr \nonumber\\
\,=&\,\vert f(0)\vert ^p\mu ([0,1))\,-\,\lim_{r\to 1}M_\infty
^p(r,f)(1-r)^{\alpha -p+1}\mu([r,1))\,+\,\int_0^1g^\prime (r)\mu
([r, 1))\,dr.\nonumber
\end{align} Since $\,f\in \Dpa\,$ we have that $\,M_\infty
^p(r,f^\prime )=\op \left ((1-r)^{-2-\alpha }\right )$, as $\,r\to
1$. Hence, $\,M_\infty ^p(r,f)=\op \left ((1-r)^{-2-\alpha +p}\right
)$, as $\,r\to 1$. This and the fact that $\,\mu \,$ is a Carleson
measure imply that
\begin{equation}\label{limit-zeroD}\lim_{r\to 1}M_\infty
^p(r,f)(1-r)^{\alpha -p+1}\mu([r,1))\,=\,0.\end{equation} Using
again that $\,\mu \,$ is a Carleson measure and integrating by parts
we see that
\begin{align*}\int_0^1g^\prime
(r)\mu ([r, 1))\,dr\,\lesssim \,&\int_0^1g^\prime (r)(1-r)\,dr\\
\,=\,&\lim_{r\to 1}g(r)(1-r)\,-g(0)\,+\,\int_0^1g(r)\,dr\\
\,\le \,& \lim_{r\to 1}M_\infty ^p(r,f)(1-r)^{\alpha -p
+2}\,+\,\int_0^1M_\infty ^p(r,f)(1-r)^{\alpha -p+1}\,dr \\
\,=\,&\int_0^1M_\infty ^p(r,f)(1-r)^{\alpha -p+1}\,dr.
\end{align*}
Then, using \cite[Lemma \,\@3]{GaGiPeSis}, it follows that
\begin{align*}\int_0^1g^\prime
(r)\mu ([r, 1))\,dr\,\lesssim \,\Vert f\Vert_{\Dpa}^p.\end{align*}
Using this and (\ref{limit-zeroD}) in (\ref{aaaaD}) readily yields
$\,\int_{[0,1)}M_\infty ^p(r,f)(1-r)^{\alpha -p+1}\,d\mu
(r)\,\lesssim \,\Vert f\Vert_{\Dpa }^p$.
\end{pf}
\par\medskip
\begin{Pf}{\,\em{Theorem\,\@\ref{Hmu-bounded-Dir}.}} Suppose first
that $\,\mathcal H_\mu \,$ is a bounded operator from $\,\Dpa \,$
into itself. For $\,1/2<b<1\,$ we set
$$f_b(z)=\frac{(1-b^2)^{1-\frac{\alpha }{p}}}{(1-bz)^{2/p}},\quad
z\in \D .$$ We have $\,\Vert f_b\Vert _{\Dpa }\asymp 1$. Then
arguing as in the proof of the correspondent implication in
Theorem\,\@\ref{Hmu-bounded-Ber} we obtain that $\,\mu \,$ is a
Carleson measure. We omit the details.
\par\medskip To prove the other implication, suppose that $\,\mu \,$
is a Carleson measure and take $\,f\in \Dpa $. Since $\,\mathcal H_
 \mu \,$ and $\,\mathcal I_\mu \,$ coincide on $\,\Dpa $, we have to
 prove that $\,\mathcal I_\mu f\in \Dpa \,$ and that $\,\Vert \mathcal I_\mu
 f\Vert_
 \Dpa \lesssim \Vert f\Vert_
 \Dpa \,$ or, equivalently, that $\,\left (\mathcal I_\mu f\right )^\prime\in A^p_\alpha \,$ and  \begin{equation}\label{Imufdpa}\,\Vert \left (
 \mathcal I_\mu
 f\right )^\prime \Vert_
 {A^p_\alpha }\lesssim \Vert f\Vert_
 {A^p_\alpha }.\end{equation}
We shall distinguish two
 cases.
\par {\bf First case: $\boldsymbol{\alpha <p-1}$.} Let $\,q\,$ be defined by the relation
$\,\frac{1}{p}+\frac{1}{q}=1\,$ and take $\,\beta =\frac{-\alpha
q}{p}$. In view of Lemma\,\@\ref{dualApa}, (\ref{Imufdpa}) is
equivalent to
\begin{equation}\label{FImufdpa}\left \vert \int_{\D }\,h(z)\,\overline{\left
(\mathcal I_\mu f\right )^\prime (z)}\,dA(z)\right \vert\,\lesssim
\Vert f\Vert _{\Dpa }\Vert h\Vert _{A^q_\beta },\quad h\in A^q_\beta
.\end{equation} So, take $\,h\in A^q_\beta$. Just as in the proof of
Theorem\,\@\ref{H1D1}, we have
\begin{equation}\label{inn}\int_{\mathbb D}\,h(z)\,\overline {\left
(\mathcal I_\mu f\right )^\prime
(z)}\,dA(z)\,=\,\int_{[0,1)}\,t\,\overline {f(t)}\,h(t)\,d\mu
(t).\end{equation} Set $\,s=-1 + \frac{\alpha+1}{p}$. Observe that
$\,ps=\alpha -p+1\,$ and $\,-qs=\beta +1$. Then, using (\ref{inn}),
H\"{o}lder's inequality, Lemma\,\@\ref{le-infty-Apa}, and
Lemma\,\@\ref{le-infty-Dpa}, we obtain
\begin{align*}&\left \vert \int_{\mathbb D}\,h(z)\,\overline {\left
(\mathcal I_\mu f\right )^\prime (z)}\,dA(z)\right\vert \,\le\,
\int_{[0,1)}\vert f(t)\vert (1-t)^s\,\vert h(t)\vert
(1-t)^{-s}\,d\mu (t)\\ \,&\le \,\left (\int_{\D }\vert f(t)\vert
^p(1-t)^{\alpha -p+1}\,d\mu (t)\right )^{1/p}\left
(\int_{[0,1)}\vert h(t)\vert ^q(1-t)^{\beta +1}\,d\mu (t)\right
)^{1/q}
\\ \,&\le \,\left (\int_{\D }M_\infty^p(t,f)(1-t)^{\alpha -p+1}\,d\mu (t)\right )^{1/p}\left
(\int_{[0,1)}M_\infty^q(t,h)(1-t)^{\beta +1}\,d\mu (t)\right )^{1/q}
\\ \,&\le \,\Vert f\Vert _{\Dpa }\Vert h\Vert_{A^q_\beta }.
\end{align*} Thus,
(\ref{FImufdpa}) holds.

\par {\bf Second case: $\boldsymbol{\alpha =p-1}$.} We let again $\,q\,$ be defined by the relation
$\,\frac{1}{p}+\frac{1}{q}=1\,$ and take $\,\beta =q-1$. Using
Lemma\,\@\ref{dualApa} and arguing as in the preceding case, we have
to show that
\begin{equation}\label{p-1q-1}\left \vert \int_{\D}\, (1-\vert
z\vert ^2)\,h(z)\,\overline {\left (\mathcal I_\mu f\right )^\prime
(z)}\,dA(z)\right \vert \,\lesssim \,\Vert f\Vert_{\Dp }\Vert h\Vert
_{A^q_{q-1}},\quad h\in A^q_{q-1}.\end{equation} We have
\begin{equation}\label{iiiii}\int_{\D}\, (1-\vert z\vert ^2)\,h(z)\,\overline {\left
(\mathcal I_\mu f\right )^\prime
(z)}\,dA(z)\,=\,\int_{[0,1)}\,t\,\overline {f(t)}\,\int_{\D
}\frac{(1-\vert z\vert^2)h(z)}{(1-t\,\overline {z})^2}\,dA(z)\,d\mu
(t).\end{equation} Now, $\int_{\D }\frac{h(z)}{(1-t\,\overline
{z})^2}\,dA(z)\,=\,h(t)$ and
\begin{align*}&\int_{\D }\frac{\vert z\vert
^2\,h(z)}{(1-t\,\overline
{z})^2}\,dA(z)\,=\,\int_0^1\,\frac{r^3}{\pi }\int_0^{2\pi
}\frac{h(re^{i\theta })\,d\theta }{(1-tre^{-i\theta })^2}\,dr\\=&\,
\int_0^1\,\frac{2r^3}{2\pi i}\int_0^{2\pi
}\frac{e^{i\theta}h(re^{i\theta })ie^{i\theta }\,d\theta
}{(e^{i\theta }-tr)^2}\,dr\,=\,\int_0^1\,\frac{2r^3}{2\pi
i}\int_{\vert z\vert
=1}\frac{zh(rz)}{(z-tr)^2}\,dz\,dr\\=&\,\int_0^12r^3\left
[h(r^2t)\,+\,r^2th^\prime (r^2t)\right ] \,dr.
\end{align*} Then it is clear that
$\left\vert \int_{\D }\frac{(1-\vert z\vert
^2)\,h(z)}{(1-t\,\overline {z})^2}\,dA(z)\right \vert \,\lesssim
M_\infty (t,h)$. Using this, (\ref{iiiii}), H\"{o}lder's inequality,
Lemma\,\@\ref{le-infty-Apa}, and Lemma\,\@\ref{le-infty-Dpa}, we
obtain
\begin{align*}&\left \vert \int_{\D}\, (1-\vert
z\vert ^2)\,h(z)\,\overline {\left (\mathcal I_\mu f\right )^\prime
(z)}\,dA(z)\right \vert \,\lesssim \,\int_{[0,1)}M_\infty
(t,f)\,M_\infty (t,h)\,d\mu(t)\\ \le \,&\left (\int_{[0,1)}M_\infty
^p(t,f)\,d\mu(t)\right )^{1/p} \,\left (\int_{[0,1)}M_\infty
^q(t,h)\,d\mu(t)\right )^{1/q} \,\le \,\Vert f\Vert_{\Dp }\Vert
h\Vert _{A^q_{q-1}}.
\end{align*} This is (\ref{p-1q-1}).
\end{Pf}
\par\medskip We shall close the article with some comments about the
case $\,p=1\,$ in Question\,\@\ref{q3}. We have the following
result.
\begin{theorem}\label{q3-p=1} Let $\,\mu \,$ be a finite positive
Borel measure on $\,[0,1)\,$ and $\,-1<\alpha <0$. If $\,\mu \,$ is
a Carleson measure then the operator $\,\mathcal H_\mu\,$ is a
bounded operator form $\,\mathcal D^1_\alpha \,$ to itself.
\end{theorem}
\begin{pf}
Using \cite[Theorem\,\@5.\,\@15, p.\,\@113]{Zhu}, we see that
 $\,A^1_\alpha \,$ can be identified as
the dual of the little Bloch space under the pairing
\begin{equation}\label{B0-a1pha}
\langle h, g\rangle\,=\,\int_{\D }(1-\vert z\vert^2)^\alpha
\,h(z)\,\overline {g(z)}\,dA(z),\quad h\in \mathcal B_0,\,\,\,g\in
A^1_\alpha .\end{equation} Suppose that $\,\mu \,$ is a Carleson
measure. Using this duality relation and the fact that $\,\mathcal
H_\mu =\mathcal I_\mu \,$ on $\,\mathcal D^1_\alpha $, showing that
$\,\mathcal H_\mu \,$ is a bounded operator from $\,\mathcal
D^1_\alpha \,$ to itself is equivalent to showing that
\begin{equation}\label{boun-D1alpha}
\left \vert \int_{\D }(1-\vert z\vert ^2)^{\alpha
}\,h(z)\,\overline{\left (\mathcal I_\mu f\right )^\prime
(z)}\,dA(z)\right\vert  \,\lesssim \Vert h\Vert _{\mathcal B}\cdot
\Vert f\Vert _{\mathcal D^1_\alpha },\quad h\in \mathcal
B_0,\,\,\,f\in \mathcal D^1_\alpha .\end{equation} Let us prove
(\ref{boun-D1alpha}). Take $\,h\in \mathcal B_0\,$ and $\,f\in
\mathcal D^1_\alpha $. We have
\begin{equation}\label{hhhhh}
\int_{\D }(1-\vert z\vert ^2)^{\alpha }\,h(z)\,\overline{\left
(\mathcal I_\mu f\right )^\prime
(z)}\,dA(z)\,=\,\int_{[0,1)}t\,\overline
{f(t)}\,\int_{\D}\frac{(1-\vert z\vert^2)^\alpha
h(z)}{(1-t\,\overline{z})^2}\,dA(z)\,d\mu (t).\end{equation} Using
\cite[Lemma\,\@5.\,\@14, pp.\,\@113-114]{Zhu} we have that the
operator $\,T\,$ defined by
$$T\phi (\xi )\,=\,(1-\vert \xi\vert ^2)^{-\alpha }\int_{\D
}\frac{(1-\vert z\vert ^2)^\alpha \phi (z)}{(1-\xi \,\overline
{z})^2}\,dA(z)$$ is a bounded operator from $\,\mathcal B\,$ into
$\,L^\infty (\D )$. Then it follows that
$$\left \vert \int_{\D}\frac{(1-\vert z\vert^2)^\alpha
h(z)}{(1-t\,\overline{z})^2}\,dA(z)\right \vert \,\lesssim \,\Vert
h\Vert _{\mathcal B}(1-t^2)^\alpha ,\quad t\in [0,1).$$ Using this
in (\ref{hhhhh}), we obtain
\begin{equation}\label{fffffff}
\left \vert \int_{\D }(1-\vert z\vert ^2)^{\alpha
}\,h(z)\,\overline{\left (\mathcal I_\mu f\right )^\prime
(z)}\,dA(z)\right\vert  \,\lesssim \Vert h\Vert _{\mathcal B}\int
_{\D }(1-t)^{\alpha }\vert f(t)\vert\,d\mu (t).\end{equation} The
fact that $\,\mu \,$ is a Carleson measure readily implies that the
measure $\,\nu \,$ defined by $\,d\nu (t)=(1-t)^\alpha \,d\mu (t)\,$
is a $\,(1-\alpha )$-Carleson measure. Using Theorem\,\@1 of
\cite{Wu} we see that then $\,\nu \,$ is a Carleson measure for
$\,\mathcal D^1_\alpha $, that is,
$$\int_{[0,1)}(1-t)^\alpha \vert g(t)\vert \,d\mu (t)\,\lesssim \,\Vert g\Vert
_{\mathcal D^1_\alpha },\quad g\in \mathcal D^1_\alpha .$$ Using
this in (\ref{fffffff}), (\ref{boun-D1alpha}) follows.
\end{pf}
\par\medskip We do not know whether the converse of
Theorem\,\@\ref{q3-p=1} is true. This is due to the fact that we do
not know whether Lemma\,\@\ref{lemDpa} remains true for $\,p=1$. The
inequality \begin{align}\label{coef-p1}\sum_{n=0}^\infty \vert
a_n\vert (n+1)^{-(1+\alpha )}\lesssim \Vert f\Vert _{\mathcal
D^1_\alpha }.\end{align}  is certainly true with no assumption on
the sequence $\,\{ a_n\} $. Indeed, by Hardy's inequality
\cite[p.\,\@48]{D}, $\sum_{n=1}^\infty \vert a_n\vert
r^{n-1}\lesssim \int_0^{2\pi }\vert f^\prime (re^{i\theta })\vert
d\theta $. Hence
\begin{align*}&\Vert f\Vert_{\mathcal D^1_\alpha }\asymp
\int_0^1(1-r)^\alpha \int_0^{2\pi }\vert f^\prime (re^{i\theta
})\vert d\theta dr \\ \gtrsim & \sum_{n=1}^\infty \vert a_n\vert
\int_0^1(1-r)^\alpha r^{n-1}\,dr = \sum_{n=1}^\infty \vert a_n\vert
B(\alpha +1, n),
\end{align*}
where $\,B(\cdot ,\cdot )\,$ is the Beta function. Stirling's
formula gives $\,B(\alpha +1, n)\asymp n^{-(\alpha +1)}$ and then
(\ref{coef-p1}) follows.
\par However, the proof of
Theorem\,\@\ref{Pav-dec-D1} in \cite{Pav-dec} does not seen to work
to prove the opposite inequality when $\{a_n\} $ is decreasing.
\par\medskip
\begin{center}{\bf Ackowledgements.}\end{center}
\par The authors wish to express their gratitude to the referees who
made several suggestions for improvement.


\end{document}